\documentclass[12pt]{amsart}

\title{Homological dimension based on a class of Gorenstein flat modules}
\author{Georgios Dalezios and Ioannis Emmanouil}
\thanks{Research supported by the Hellenic Foundation for 
Research and Innovation (H.F.R.I.) under the ``1st Call 
for H.F.R.I. Research Projects to support Faculty members 
and Researchers and the procurement of high-cost research 
equipment grant", project number 4226.}

\newtheorem{Lemma}{Lemma}[section]
\newtheorem{Proposition}[Lemma]{Proposition}
\newtheorem{Theorem}[Lemma]{Theorem}
\newtheorem{Corollary}[Lemma]{Corollary}

\newcommand{\Limit}{\mbox{$\displaystyle{\lim_{\longleftarrow}}$}}

\begin{document}

\begin{abstract}
In this paper, we study the relative homological dimension based 
on the class of projectively coresolved Gorenstein flat modules 
(PGF-modules), that were introduced by Saroch and Stovicek in [26]. 
The resulting PGF-dimension of modules has several properties in 
common with the Gorenstein projective dimension, the relative 
homological theory based on the class of Gorenstein projective 
modules. In particular, there is a hereditary Hovey triple in the 
category of modules of finite PGF-dimension, whose associated 
homotopy category is triangulated equivalent to the stable 
category of PGF-modules. Studying the finiteness of the PGF 
global dimension reveals a connection between classical homological 
invariants of left and right modules over the ring, that leads 
to generalizations of certain results by Jensen [24], Gedrich 
and Gruenberg [17] that were originally proved in the realm of 
commutative Noetherian rings.
\end{abstract}

\maketitle
\tableofcontents

\addtocounter{section}{-1}
\section{Introduction}

\noindent
The concept of G-dimension for commutative Noetherian rings that was 
introduced by Auslander and Bridger in [1] has been extended to modules 
over any ring $R$ through the notion of a Gorenstein projective module. 
Such a module is, by definition, a syzygy of an acyclic complex of 
projective modules which remains acyclic when applying the functor 
$\mbox{Hom}_{R}(\_\!\_ ,P)$ for any projective module $P$. The modules 
of finite Gorenstein projective dimension are defined in the standard way, 
using resolutions by Gorenstein projective modules. A Gorenstein flat
module is a syzygy of an acyclic complex of flat modules which remains 
acyclic when applying the functor $I \otimes_R \_\!\_$ for any injective
right module $I$. The modules of finite Gorenstein flat dimension are 
then defined using resolutions by Gorenstein flat modules. The standard
reference for these notions is Holm's paper [22]. The relation between 
Gorenstein projective and Gorenstein flat modules remains somehow
mysterious in general. As shown in [loc.cit.], all Gorenstein projective 
modules are Gorenstein flat if the ground ring is right coherent and has 
finite left finitistic dimension (i.e.\ if there is an upper bound on the 
projective dimension of all modules that have finite projective dimension).

The projectively coresolved Gorenstein flat modules (PGF-modules, for short) 
were defined by Saroch and Stovicek in [26]; these are the syzygies of the 
acyclic complexes of projective modules that remain acyclic when applying 
the functor $I \otimes_R \_\!\_$ for any injective right module $I$. It is 
clear that PGF-modules are Gorenstein flat. As shown in [26, Theorem 4.4],
PGF-modules are also Gorenstein projective. A schematic presentation of the 
classes ${\tt GProj}(R)$, ${\tt GFlat}(R)$ and ${\tt PGF}(R)$ of Gorenstein 
projective, Gorenstein flat and PGF-modules respectively is given below
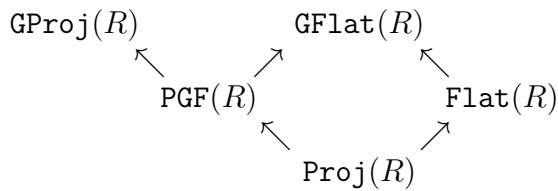
\[
\begin{array}{ccccccc}
 {\tt GProj}(R) \!\!\! & & & & \!\!\! {\tt GFlat}(R) \!\!\! & & \\
 & \!\!\! \nwarrow \!\!\! & & \!\!\! \nearrow& \!\!\! & \!\!\! 
 \nwarrow \!\!\! & \\
 & & \!\!\! {\tt PGF}(R) \!\!\! & & & & \!\!\! {\tt Flat}(R) \!\!\! \\
 & & & \!\!\! \nwarrow \!\!\! & & \!\!\! \nearrow \!\!\! & \\
 & & & & \!\!\! {\tt Proj}(R) \!\!\! & & 
\end{array}
\]
Here, ${\tt Proj}(R)$ and ${\tt Flat}(R)$ denote the classes of projective 
and flat modules respectively and all arrows are inclusions. Moreover, the 
class ${\tt Proj}(R)$ of projective modules is the intersection 
${\tt PGF}(R) \cap {\tt Flat}(R)$ and all classes pictured above are 
projectively resolving; in fact, ${\tt GFlat}(R)$ is the smallest 
projectively resolving class of modules that contains both 
${\tt PGF}(R)$ and ${\tt Flat}(R)$; these assertions are proved in [26].

In this paper, we study the relative homological dimension which is based 
on the class ${\tt PGF}(R)$ and define the PGF-dimension $\mbox{PGF-dim}_RM$ 
of a module $M$ as the minimal length of a resolution of $M$ by PGF-modules 
(provided that such a resolution exists). The resulting class 
$\overline{\tt PGF}(R)$ of modules of finite PGF-dimension has many of the 
standard properties that one would expect. In particular, it is closed under 
direct sums, direct summands and has the 2-out-of-3 property for short exact 
sequences of modules. The PGF-dimension is a refinement of the ordinary projective 
dimension, whereas the Gorenstein projective dimension is a refinement of the 
PGF-dimension. In other words, if $M$ is a module of finite projective dimension 
(resp.\ of finite PGF-dimension), then $M$ has finite PGF-dimension (resp.\ 
finite Gorenstein projective dimension) and $\mbox{PGF-dim}_RM = \mbox{pd}_RM$ 
(resp.\ $\mbox{Gpd}_RM = \mbox{PGF-dim}_RM$). When restricted to the class 
${\tt Flat}(R)$ of flat modules, the PGF-dimension coincides with the projective 
dimension. The modules of finite PGF-dimension can be approximated by modules 
of finite projective dimension and PGF-modules, in analogy with the case of 
modules of finite Gorenstein projective dimension. In particular, this leads 
to a description, up to triangulated equivalence, of the stable category of 
PGF-modules modulo projective modules, as the homotopy category of the exact 
model structure which is associated with a Hovey triple in the category 
$\overline{\tt PGF}(R)$. Using the analogous approximations of Gorenstein 
flat modules by PGF-modules and flat modules, that were obtained by Saroch 
and Stovicek in [26], we describe a similar Hovey triple in the category 
${\tt GFlat}(R)$. Therefore, in order to realize the stable category of 
PGF-modules as the homotopy category of a Quillen model structure, it is
sufficient to work on either subcategory $\overline{\tt PGF}(R)$ or 
${\tt GFlat}(R)$ of the module category.

In order to present an application of the notion of PGF-dimension studied
in this paper, we consider the invariants $\mbox{silp} \, R$ and 
$\mbox{spli} \, R$, which are defined as the suprema of the injective 
lengths of projective modules and the projective lengths of injective 
modules, respectively. It is easily seen that these invariants are equal, 
if they are both finite. Nevertheless, as Gedrich and Gruenberg point out 
in [17], it is not clear whether the finiteness of one of these implies the
finiteness of the other, i.e.\ whether we always have an equality
$\mbox{silp} \, R = \mbox{spli} \, R$. In the special case where $R$ is an
Artin algebra, the equality $\mbox{silp} \, R = \mbox{spli} \, R$ is
equivalent to the Gorenstein Symmetry Conjecture in representation theory;
cf.\ [2, Conjecture 13], [4, $\S 11$] and [5, Chapter VII]. 

The study of the finiteness of the PGF global dimension reveals a connection 
between the silp and spli invariants for left and right modules over any ring, 
which may be itself used in order to show that:
\begin{center}
 {\em If both $\mbox{spli} \, R$ and $\mbox{spli} \, R^{op}$ 
 are finite, then $\mbox{silp} \, R = \mbox{spli} \, R$ and 
 $\mbox{silp} \, R^{op} = \mbox{spli} \, R^{op}$.}
\end{center}
Using the Hopf algebra structure of the group algebra $kG$ of a group $G$ 
with coefficients in a commutative ring $k$, Gedrich and Gruenberg proved
in [17] that $\mbox{silp} \, kG \leq \mbox{spli} \, kG$, in the special 
case where the commutative ring $k$ is Noetherian of finite self-injective
dimension. It follows from the 
result displayed above that we actually have an inequality 
$\mbox{silp} \, R \leq \mbox{spli} \, R$ for any ring $R$ which is 
isomorphic with its opposite ring $R^{op}$. In particular, the inequality 
holds for group algebras of groups over {\em any} commutative coefficient 
ring. On the other hand, Jensen has proved in [24, 5.9] that the equality 
$\mbox{silp} \, R = \mbox{spli} \, R$ holds for any commutative Noetherian 
ring $R$. The result displayed above, combined with earlier work in [14], 
shows that the equality $\mbox{silp} \, R = \mbox{spli} \, R$ actually 
holds for any commutative $\aleph_0$-Noetherian ring $R$, i.e.\ for any 
commutative ring $R$ all of whose ideals are countably generated.
\vspace{0.1in}
\newline
{\em Notations and terminology.}
We work over a fixed unital associative ring $R$ and, unless otherwise 
specified, all modules are left $R$-modules. We denote by $R^{op}$ the 
opposite ring of $R$ and do not distinguish between right $R$-modules 
and left $R^{op}$-modules. If $\lambda(R)$ is an invariant, which is 
defined in terms of a certain class of left $R$-modules, then we denote 
by $\lambda(R^{op})$ the corresponding invariant, which is defined for 
$R$ in terms of the appropriate class of right $R$-modules. Finally, we 
say that a class ${\mathcal C}$ of modules is projectively resolving if
${\tt Proj}(R) \subseteq {\mathcal C}$ and ${\mathcal C}$ is closed under 
extensions and kernels of epimorphisms.

\section{Preliminary notions}

\noindent
In this section, we collect certain basic notions and preliminary 
results that will be used in the sequel. These involve basic concepts
related to Gorenstein homological algebra in module categories and the
theory of Hovey triples in exact additive categories.
\vspace{0.15in}
\newline
{\sc I.\ Gorenstein projective and Gorenstein flat modules.}
An acyclic complex $P_*$ of projective modules is said to be a complete 
projective resolution if the complex of abelian groups $\mbox{Hom}_R(P_*,Q)$
is acyclic for any projective module $Q$. Then, a module is called Gorenstein 
projective if it is a syzygy of a complete projective resolution. Holm's 
paper [22] is the standard reference in Gorenstein homological algebra. 
The class ${\tt GProj}(R)$ of Gorenstein projective modules is projectively 
resolving; it is also closed under direct sums and direct summands. The 
Gorenstein projective dimension $\mbox{Gpd}_RM$ of a module $M$ is the 
length of a shortest resolution of $M$ by Gorenstein projective modules. 
If no such resolution of finite length exists, then we write 
$\mbox{Gpd}_RM = \infty$. If $M$ is a module of finite projective 
dimension, then $M$ has finite Gorenstein projective dimension as 
well and $\mbox{Gpd}_RM = \mbox{pd}_RM$.

An acyclic complex $F_*$ of flat modules is said to be a complete flat
resolution if the complex of abelian groups $I \otimes_R F_*$ is acyclic 
for any injective right module $I$. We say that a module is Gorenstein 
flat if it is a syzygy of a complete flat resolution. We let ${\tt GFlat}(R)$ 
be the class of Gorenstein flat modules. The Gorenstein flat dimension 
$\mbox{Gfd}_RM$ of a module $M$ is the length of a shortest resolution 
of $M$ by Gorenstein flat modules. If no such resolution of finite length
exists, then we write $\mbox{Gfd}_RM = \infty$. If $M$ is a module of 
finite flat dimension, then $M$ has finite Gorenstein flat dimension as
well and $\mbox{Gfd}_RM = \mbox{fd}_RM$.

Even though the relation between Gorenstein projective and Gorenstein 
flat modules is not fully understood, the notion of a projectively 
coresolved Gorenstein flat module (for short, PGF-module) defined in 
[26] sheds some light in and helps clarifying that relation. A PGF-module 
is a syzygy of an acyclic complex of projective modules $P_*$, which is 
such that the complex of abelian groups $I \otimes_R P_*$ is acyclic for 
any injective right module $I$. It is clear that the class ${\tt PGF}(R)$ 
of PGF-modules is contained in ${\tt GFlat}(R)$. The inclusion 
${\tt PGF}(R) \subseteq {\tt GProj}(R)$ is proved in [26, Theorem 4.4]; 
in fact, it is shown that $\mbox{Ext}^1_R(M,F)=0$ for any PGF-module 
$M$ and any flat module $F$. It is also proved in [loc.cit.] that the 
classes ${\tt PGF}(R)$ and ${\tt GFlat}(R)$ are both projectively 
resolving, closed under direct sums and direct summands. 
\vspace{0.15in}
\newline
{\sc II.\ Gorenstein global dimensions.}
The existence of complete projective resolutions of modules (i.e.\ of
complete projective resolutions that coincide in sufficiently large 
degrees with an ordinary projective resolution of the module) has been 
studied by Gedrich and Gruenberg [17], Cornick and Kropholler [12], in 
connection with the existence of complete cohomological functors in the 
category of modules. Even though they were mainly interested in the case 
where $R$ is the integral group ring of a group, they were able to 
characterize those rings over which all modules admit complete projective 
resolutions, in terms of the finiteness of the invariants $\mbox{spli} \, R$ 
and $\mbox{silp} \, R$. Here, $\mbox{spli} \, R$ is the supremum of the 
projective lengths (dimensions) of injective modules and $\mbox{silp} \, R$ 
is the supremum of the injective lengths (dimensions) of projective modules. 
As shown by Holm [22], the existence of a complete projective resolution for 
a module $M$ is equivalent to the finiteness of the Gorenstein projective 
dimension $\mbox{Gpd}_RM$ of $M$. From this point of view, the above result 
by Cornick and Kropholler was alternatively proved by Bennis and Mahbou in 
[7], where the notion of the Gorenstein global dimension of the ring was 
introduced, in analogy with the classical notion of global dimension defined 
in [9, Chapter VI, $\S $2]; see also [15, $\S $4]. More precisely, the (left) 
Gorenstein global dimension $\mbox{Ggl.dim} \, R$ of the ring $R$ is defined 
by letting 
\[ \mbox{Ggl.dim} \, R = \sup \{ \mbox{Gpd}_{R}M : 
   \mbox{$M$ a left $R$-module} \} . \]
Then, the following conditions are equivalent:

(i) $\mbox{Ggl.dim} \, R < \infty$, 

(ii) $\mbox{Gpd}_{R}M < \infty$ for any module $M$,

(iii) any module $M$ admits a complete projective resolution,

(iv) the invariants $\mbox{spli} \, R$ and $\mbox{silp} \, R$ are
     finite.
\newline
If these conditions are satisfied, then
$\mbox{Ggl.dim} \, R =  \mbox{spli} \, R = \mbox{silp} \, R$.

The corresponding characterization of the finiteness of the 
Gorenstein weak global dimension $\mbox{Gwgl.dim} \, R$ of 
the ring $R$, which is defined by letting 
\[ \mbox{Gwgl.dim} \, R = \sup \{ \mbox{Gfd}_{R}M : 
   \mbox{$M$ a left $R$-module} \} , \]
turned out to be more difficult to achieve. The relevant homological
invariants here are $\mbox{sfli} \, R$, the supremum of the flat 
lengths (dimensions) of injective modules, and its analogue 
$\mbox{sfli} \, R^{op}$ for the opposite ring $R^{op}$. Using in 
an essential way results in [26], it was proved by Christensen, 
Estrada and Thompson in [11] that the following conditions are equivalent:

(i) $\mbox{Gwgl.dim} \, R < \infty$, 

(ii) $\mbox{Gfd}_{R}M < \infty$ for any module $M$,

(iii) the invariants $\mbox{sfli} \, R$ and $\mbox{sfli} \, R^{op}$ 
      are finite.
\newline
If these conditions are satisfied, then
$\mbox{Gwgl.dim} \, R =  \mbox{sfli} \, R = \mbox{sfli} \, R^{op}$.
\vspace{0.15in}
\newline
{\sc III.\ Cotorsion pairs and Hovey triples.}
Let ${\mathcal A}$ be an exact additive category, in the sense of 
Quillen [8], and consider a full subcategory 
${\mathcal B} \subseteq {\mathcal A}$. A morphism 
$f : B \longrightarrow A$ in ${\mathcal A}$ is called a 
${\mathcal B}$-precover of the object $A \in {\mathcal A}$ if:

(i) $B \in {\mathcal B}$ and 

(ii) the induced map 
     $f_* : \mbox{Hom}_{\mathcal A}(B',B) \longrightarrow 
            \mbox{Hom}_{\mathcal A}(B',A)$
     is surjective for any $B' \in {\mathcal B}$.
\newline
The reader is referred to [20] for a thorough and systematic study 
of precovers.

The $\mbox{Ext}^1$-pairing induces an orthogonality relation between 
subclasses of ${\mathcal A}$. If ${\mathcal B} \subseteq {\mathcal A}$, 
then we define the left orthogonal $^{\perp}{\mathcal B}$ of 
${\mathcal B}$ as the class consisting of those objects 
$X \in {\mathcal A}$, which are such that 
$\mbox{Ext}^1_{\mathcal A}(X,B)=0$ for all $B \in {\mathcal B}$. 
Analogously, the right orthogonal ${\mathcal B}^{\perp}$ of 
${\mathcal B}$ is the class consisting of those objects 
$Y \in {\mathcal A}$, which are such that 
$\mbox{Ext}^1_{\mathcal A}(B,Y)=0$ for all $B \in {\mathcal B}$. 
If ${\mathcal C},{\mathcal D}$ are two subclasses of ${\mathcal A}$, 
then the pair $({\mathcal C},{\mathcal D})$ is a cotorsion pair in
${\mathcal A}$ (cf.\ [16]) if ${\mathcal C} = {^{\perp} {\mathcal D}}$ 
and ${\mathcal C} ^{\perp} = {\mathcal D}$. The cotorsion pair is 
called hereditary if $\mbox{Ext}^i_{\mathcal A}(C,D)=0$ for all $i>0$
and all objects $C \in {\mathcal C}$ and $D \in {\mathcal D}$. The 
cotorsion pair is complete if for any object $A \in {\mathcal A}$ 
there exist short exact sequences (conflations), usually called 
approximation sequences
\[ 0 \longrightarrow D \longrightarrow C \longrightarrow A 
     \longrightarrow 0 
   \;\;\; \mbox{and } \;\;\;
   0 \longrightarrow A \longrightarrow D' \longrightarrow C' 
     \longrightarrow 0 , \]
where $C,C' \in {\mathcal C}$ and $D,D' \in {\mathcal D}$. In that 
case, the morphism $C \longrightarrow A$ is a ${\mathcal C}$-precover 
of $A$.

A Hovey triple on ${\mathcal A}$ is a triple 
$({\mathcal C} , {\mathcal W} , {\mathcal F})$ of subclasses of 
${\mathcal A}$, which are such that the pairs 
$({\mathcal C} , {\mathcal W} \cap {\mathcal F})$ and
$({\mathcal C} \cap {\mathcal W} , {\mathcal F})$ are complete 
cotorsion pairs and the class ${\mathcal W}$ is closed under direct 
summands and satisfies the 2-out-of-3 property for short exact 
sequences (conflations) in ${\mathcal A}$. The fundamental work 
of Gillespie [18], which is based on work of Hovey [23], gives 
a bijection between Hovey triples on a (weakly) idempotent complete 
exact category ${\mathcal A}$ and certain, so-called exact, Quillen 
model structures on ${\mathcal A}$; cf.\ [18, Theorem 3.3]. In the 
context of Gillespie's bijection, it is proved in [18, Proposition 
5.2] that for an exact model structure on ${\mathcal A}$ that has 
its associated complete cotorsion pairs hereditary, the class 
${\mathcal C} \cap {\mathcal F}$ is a Frobenius exact category with 
projective-injective objects equal to 
${\mathcal C} \cap {\mathcal W} \cap {\mathcal F}$. Then, a result 
of Happel [21] implies that the associated stable category, which 
is ${\mathcal C} \cap {\mathcal F}$ modulo its projective-injective 
objects, is triangulated. The upshot of this connection is that the 
(Quillen) homotopy category of an exact model structure is triangulated 
equivalent to the stable category of the Frobenius exact category 
${\mathcal C} \cap {\mathcal F}$; cf.\ [18, Proposition 4.4 and 
Corollary 4.8].

\section{Modules of finite PGF-dimension}

\noindent
In this section, we define the notion of PGF-dimension for a module 
and show that the resulting class $\overline{\tt PGF}(R)$ of modules 
of finite PGF-dimension has many standard closure properties.

We recall that the class ${\tt PGF}(R)$ is projectively resolving and 
closed under direct sums and direct summands. The following result is 
a formal consequence of these properties of ${\tt PGF}(R)$; cf.\ 
[1, Lemma 3.12].

\begin{Lemma}
Let $M$ be an $R$-module, $n$ a non-negative integer and 
\[ 0 \longrightarrow K \longrightarrow G_{n-1} \longrightarrow \cdots
     \longrightarrow G_0 \longrightarrow M \longrightarrow 0 , \]
\[ 0 \longrightarrow K' \longrightarrow G'_{n-1} \longrightarrow \cdots
     \longrightarrow G'_0 \longrightarrow M \longrightarrow 0 \]
two exact sequences of modules with 
$G_0, \ldots , G_{n-1},G'_0, \ldots , G'_{n-1} \in {\tt PGF}(R)$.
Then, $K \in {\tt PGF}(R)$ if and only if $K' \in {\tt PGF}(R)$.
\hfill $\Box$
\end{Lemma}

\begin{Proposition}
The following conditions are equivalent for an $R$-module $M$ and 
a non-negative integer $n$:

(i) There exists an exact sequence of modules
\[ 0 \longrightarrow G_n \longrightarrow G_{n-1} \longrightarrow \cdots
     \longrightarrow G_0 \longrightarrow M \longrightarrow 0 , \]
with $G_0, \ldots , G_{n-1},G_n \in {\tt PGF}(R)$.

(ii) For any exact sequence of modules
\[ 0 \longrightarrow K \longrightarrow G_{n-1} \longrightarrow \cdots
     \longrightarrow G_0 \longrightarrow M \longrightarrow 0 \]
with $G_0, \ldots , G_{n-1} \in {\tt PGF}(R)$, we also have 
$K \in {\tt PGF}(R)$.
\end{Proposition}
\vspace{-0.05in}
\noindent 
{\em Proof.}
The implication (i)$\rightarrow$(ii) is a consequence of Lemma 2.1, 
whereas the implication (ii)$\rightarrow(i)$ follows by considering 
a truncated projective resolution of $M$. \hfill $\Box$ 
\vspace{0.15in}
\newline
If the equivalent conditions in Proposition 2.2 are satisfied, then 
we say that the module $M$ has a PGF-resolution of length $n$ and 
write $\mbox{PGF-dim}_RM \leq n$. In the case where 
$\mbox{PGF-dim}_RM \leq n$ and $M$ has no PGF-resolution of length 
$<n$, we say that $M$ has PGF-dimension equal to $n$ and write 
$\mbox{PGF-dim}_RM=n$. Finally, we say that $M$ has infinite PGF-dimension 
and write $\mbox{PGF-dim}_RM = \infty$, if $M$ has no PGF-resolution 
of finite length. 

We now consider the class $\overline{\tt PGF}(R)$ of all modules of 
finite PGF-dimension and describe certain closure properties of that 
class.

\begin{Proposition}
Let $(M_i)_i$ be a family of modules and $M = \bigoplus_iM_i$ the 
corresponding direct sum. Then, 
$\mbox{PGF-dim}_RM = \sup_i \mbox{PGF-dim}_RM_i$. In particular, 
the class $\overline{\tt PGF}(R)$ is closed under finite direct 
sums and direct summands.
\end{Proposition}
\vspace{-0.05in}
\noindent
{\em Proof.}
In order to show that 
$\mbox{PGF-dim}_RM \leq \sup_i \mbox{PGF-dim}_RM_i$, it suffices to 
consider the case where $\sup_i \mbox{PGF-dim}_RM_i = n < \infty$. 
Then, $\mbox{PGF-dim}_RM_i \leq n$ and hence $M_i$ has a PGF-resolution 
of length $n$ for all $i$. Since ${\tt PGF}(R)$ is closed under direct
sums, the direct sum of these resolutions is a PGF-resolution of $M$ of 
length $n$, so that $\mbox{PGF-dim}_RM \leq n$.

It remains to show that we also have
$\sup_i \mbox{PGF-dim}_RM_i \leq \mbox{PGF-dim}_RM$. To that end, 
assume that $\mbox{PGF-dim}_RM = n < \infty$ and consider for any 
$i$ an exact sequence 
\[ 0 \longrightarrow K_i \longrightarrow G_{i,n-1} \longrightarrow \cdots
     \longrightarrow G_{i,0} \longrightarrow M_i \longrightarrow 0 , \]
with $G_{i,0}, \ldots , G_{i,n-1} \in {\tt PGF}(R)$. Since ${\tt PGF}(R)$
is closed under direct sums, the exactness of the direct sum of these 
exact sequences
\[ 0 \longrightarrow {\textstyle{\bigoplus_i}} K_i 
     \longrightarrow {\textstyle{\bigoplus_i}} G_{i,n-1} 
     \longrightarrow \cdots 
     \longrightarrow {\textstyle{\bigoplus_i}} G_{i,0} 
     \longrightarrow M \longrightarrow 0 \]
and our assumption on the PGF-dimension of $M$ imply that 
${\textstyle{\bigoplus_i}} K_i$ is a PGF-module. Since ${\tt PGF}(R)$
is closed under direct summands, it follows that $K_i$ is a PGF-module 
for all $i$. Then, $M_i$ has a PGF-resolution of length $n$ and hence 
$\mbox{PGF-dim}_RM_i \leq n$ for all $i$, as needed. \hfill $\Box$

\begin{Proposition}
Let  
$0 \longrightarrow M' \longrightarrow M \longrightarrow M'' 
   \longrightarrow 0$
be a short exact sequence of modules. Then:

(i) 
$\mbox{PGF-dim}_RM \leq 
 \max \{ \mbox{PGF-dim}_RM' , \mbox{PGF-dim}_RM'' \}$,
 
(ii) 
$\mbox{PGF-dim}_RM' \leq 
 \max \{ \mbox{PGF-dim}_RM , \mbox{PGF-dim}_RM'' \}$,

(iii) 
$\mbox{PGF-dim}_RM'' \leq 1 +  
 \max \{ \mbox{PGF-dim}_RM' , \mbox{PGF-dim}_RM \}$.
\newline
In particular, the class $\overline{\tt PGF}(R)$ has the 2-out-of-3 
property: if two out of the three modules that appear in a short exact 
sequence have finite PGF-dimension, then so does the third.
\end{Proposition}
\vspace{-0.05in}
\noindent
{\em Proof.}
(i) Assume that 
$\max \{ \mbox{PGF-dim}_RM' , \mbox{PGF-dim}_RM'' \} = n$
and consider two projective resolutions 
$P'_* \longrightarrow M' \longrightarrow 0$ and 
$P''_* \longrightarrow M'' \longrightarrow 0$ of $M'$ and
$M''$ respectively. Then, we may construct by the standard 
step-by-step process a projective resolution 
$P_* \longrightarrow M \longrightarrow 0$ of $M$, such that
$P_i = P'_i \oplus P''_i$ and the corresponding syzygy module
$\Omega_iM$ is an extension of $\Omega_iM''$ by $\Omega_iM'$
for all $i$. Since both $M'$ and $M''$ have PGF-dimension 
$\leq n$, the modules $\Omega_nM'$ and $\Omega_nM''$ are both 
PGF-modules. Then, the short exact sequence 
\[ 0 \longrightarrow \Omega_nM' \longrightarrow \Omega_nM 
     \longrightarrow \Omega_nM'' \longrightarrow 0 \]
and the closure of ${\tt PGF}(R)$ under extensions show that 
$\Omega_nM$ is a PGF-module as well. Then, the exact sequence
\[ 0 \longrightarrow \Omega_nM \longrightarrow P_{n-1} 
     \longrightarrow \cdots \longrightarrow P_o
     \longrightarrow M \longrightarrow 0 \]
is a PGF-resolution of $M$ of length $n$ and hence 
$\mbox{PGF-dim}_RM \leq n$, as needed.

(ii) We can prove this assertion by using the same argument as 
the one used in order to prove assertion (i) above, by invoking 
the closure of ${\tt PGF}(R)$ under kernels of epimorphisms.

(iii) We fix a short exact sequence 
\begin{equation}
 0 \longrightarrow K \longrightarrow P
   \stackrel{p}{\longrightarrow} M'' \longrightarrow 0 , 
\end{equation}
where $P$ is a projective module, and consider the pullback 
of the short exact sequence given in the statement of the 
Proposition along $p$
\[
\begin{array}{ccccccccc}
 & & & & 0 & & 0 & \\
 & & & & \downarrow & & \downarrow & \\
 & & & & K & = & K & \\
 & & & & \downarrow & & \downarrow & \\
 0 & \longrightarrow & M' & \longrightarrow & X 
   & \longrightarrow & P & \longrightarrow & 0 \\
 & & \parallel & & \downarrow 
 & & \!\!\!\! {\scriptstyle{p}} \downarrow & \\
 0 & \longrightarrow & M' & \longrightarrow & M 
   & \longrightarrow & M'' & \longrightarrow & 0 \\
 & & & & \downarrow & & \downarrow & \\
 & & & & 0 & & 0 &
\end{array}
\]
Since $P$ is projective, the horizontal short exact sequence 
in the middle of the diagram splits and hence $X \simeq P \oplus M'$. 
We now invoke Proposition 2.3 and conclude that 
$\mbox{PGF-dim}_RX = \mbox{PGF-dim}_RM'$. Then, the vertical 
short exact sequence in the middle of the diagram and assertion 
(ii) above show that
\[ \mbox{PGF-dim}_RK \leq 
   \max \{ \mbox{PGF-dim}_RX , \mbox{PGF-dim}_RM \} =
   \max \{ \mbox{PGF-dim}_RM' , \mbox{PGF-dim}_RM \} . \]
Since we may splice any PGF-resolution of $K$ of length 
$\mbox{PGF-dim}_RK$ with the short exact sequence (1) and obtain 
a PGF-resolution of $M''$ of length $1+\mbox{PGF-dim}_RK$, it 
follows that 
\[ \mbox{PGF-dim}_RM'' \leq 1 + \mbox{PGF-dim}_RK \leq   
   1+ \max \{ \mbox{PGF-dim}_RM' , \mbox{PGF-dim}_RM \} , \]
as needed. \hfill $\Box$
\vspace{0.15in}
\newline
As a consequence of the equality  
${\tt PGF}(R) \cap {\tt Flat}(R) = {\tt Proj}(R)$, we obtain the 
following result on the relation between the projective dimension 
and the PGF-dimension of flat modules and, analogously, the relation 
between the projective dimension and the flat dimension of PGF-modules.

\begin{Proposition}
(i) If $M$ is a flat module, then $\mbox{pd}_RM = \mbox{PGF-dim}_RM$.

(ii) If $M$ is a PGF-module, then $\mbox{pd}_RM = fd_RM$.
\end{Proposition}
\vspace{-0.05in}
\noindent
{\em Proof.}
(i) Since ${\tt Proj}(R) \subseteq {\tt PGF}(R)$, we always have 
$\mbox{PGF-dim}_RM \leq \mbox{pd}_RM$. In order to prove the reverse 
inequality, it suffices to assume that $\mbox{PGF-dim}_RM = n < \infty$.
Then, the truncation of a projective resolution of $M$ provides us with 
an exact sequence
\[ 0 \longrightarrow K \longrightarrow P_{n-1} 
     \longrightarrow \cdots \longrightarrow P_0 
     \longrightarrow M \longrightarrow 0 , \]
where $P_0, \ldots ,P_{n-1}$ are projective modules and 
$K \in {\tt PGF}(R)$. Since $M$ is flat, it follows that $K$ is also 
flat and hence $K \in {\tt PGF}(R) \cap {\tt Flat}(R) = {\tt Proj}(R)$.
We conclude that $M$ admits a projective resolution of length $n$ and 
hence $\mbox{pd}_RM \leq n = \mbox{PGF-dim}_RM$, as needed.

(ii) Since projective modules are flat, we always have 
$\mbox{fd}_RM \leq \mbox{pd}_RM$. In order to prove the reverse 
inequality, it suffices to assume that $\mbox{fd}_RM = n < \infty$. 
Then, the truncation of a projective resolution of $M$ provides us 
with an exact sequence
\[ 0 \longrightarrow K \longrightarrow P_{n-1} 
     \longrightarrow \cdots \longrightarrow P_0 
     \longrightarrow M \longrightarrow 0 , \]
where $P_0, \ldots , P_{n-1}$ are projective modules and $K$ is flat. 
Since $M$ is a PGF-module and the class ${\tt PGF}(R)$ is projectively 
resolving, it follows that $K$ is also a PGF-module. Then, 
$K \in {\tt PGF}(R) \cap {\tt Flat}(R) = {\tt Proj}(R)$ and hence $M$ 
admits a projective resolution of length $n$, i.e.\ 
$\mbox{pd}_RM \leq n = \mbox{fd}_RM$. \hfill $\Box$
\vspace{0.15in}
\newline
{\bf Remark 2.6.}
If we denote by $\overline{\tt Proj}(R)$ and $\overline{\tt Flat}(R)$ 
the classes of modules of finite projective dimension and finite flat 
dimension respectively, then 
$\overline{\tt PGF}(R) \cap \overline{\tt Flat}(R) =
 \overline{\tt Proj}(R)$.
Indeed, the inclusion 
$\overline{\tt Proj}(R) \subseteq 
 \overline{\tt PGF}(R) \cap \overline{\tt Flat}(R)$
is clear, since any projective resolution of finite length is both 
a PGF-resolution and a flat resolution of finite length. Conversely, 
if $M$ is a module contained in 
$\overline{\tt PGF}(R) \cap \overline{\tt Flat}(R)$, then the $n$-th
syzygy module $\Omega_nM$ in a projective resolution of $M$ is a flat 
and PGF-module for $n \gg 0$. Since 
${\tt PGF}(R) \cap {\tt Flat}(R) = {\tt Proj}(R)$, it follows that 
$\Omega_nM$ is projective for $n \gg 0$ and hence 
$M \in \overline{\tt Proj}(R)$.
\addtocounter{Lemma}{1}

\section{Approximation sequences}

\noindent
In this section, we show that the finiteness of PGF-dimension can 
be detected by the existence of suitable approximation sequences, 
in analogy with the case of the finiteness of Gorenstein projective 
dimension.

The next result is akin to [22, Theorem 2.10].

\begin{Proposition}
Let $M$ be a module with $\mbox{PGF-dim}_RM = n$. Then, there exists 
a short exact sequence 
\[ 0 \longrightarrow K \longrightarrow G 
     \stackrel{\pi}{\longrightarrow} M \longrightarrow 0 , \]
where $G$ is a PGF-module and $\mbox{pd}_RK = n-1$. (If $n=0$, 
this is understood to mean $K=0$.) In particular, $\pi$ is a 
${\tt PGF}(R)$-precover of $M$.  
\end{Proposition}
\vspace{-0.05in}
\noindent
{\em Proof.}
The result is clear if $n=0$ and hence we may assume that $n \geq 1$.
Since $\mbox{PGF-dim}_RM = n$, there exists an exact sequence
\[ 0 \longrightarrow N \longrightarrow P_{n-1} \longrightarrow 
     \cdots \longrightarrow P_0 \longrightarrow M \longrightarrow 0 \]
where $P_0, \ldots ,P_{n-1}$ are projective modules and $N \in {\tt PGF}(R)$.
Then, there exists another exact sequence 
\[ 0 \longrightarrow N \longrightarrow Q_0 \longrightarrow Q_{-1}
     \longrightarrow \cdots \longrightarrow Q_{-n+1} 
     \longrightarrow G \longrightarrow 0 , \]
where $Q_0, \ldots ,Q_{-n+1}$ are projective modules and 
$G \in {\tt PGF}(R)$. 
Since all kernels of the latter exact sequence are PGF-modules as 
well, it follows from [26, Corollary 4.5] that the exact sequence 
remains exact after applying the functor $\mbox{Hom}_R(\_\!\_,P)$ 
for any projective module $P$. We conclude that there exists a 
morphism of complexes
\[
\begin{array}{ccccccccccc}
0 & \longrightarrow & N & \longrightarrow & Q_0 
  & \longrightarrow \cdots \longrightarrow & Q_{-n+1} 
  & \longrightarrow & G & \longrightarrow & 0 \\
  & & \parallel & & \downarrow & & \downarrow 
  & & \downarrow & & \\ 
0 & \longrightarrow & N & \longrightarrow & P_{n-1} 
  & \longrightarrow \cdots \longrightarrow & P_0 
  & \longrightarrow & M & \longrightarrow & 0
\end{array}
\]
The unlabelled vertical arrows induce a quasi-isomorphism between 
the corresponding complexes and hence we may consider the 
associated mapping cone, which is an acyclic complex
\begin{equation}
 0 \longrightarrow Q_0 \longrightarrow Q_{-1} \oplus P_{n-1}
   \longrightarrow \cdots \longrightarrow G \oplus P_0
   \stackrel{\pi}{\longrightarrow} M \longrightarrow 0 . 
\end{equation}
Note that $G \oplus P_0$ is a PGF-module and the module 
$K = \mbox{ker} \, \pi$ has projective dimension $\leq n-1$. 
In fact, our assumption that $\mbox{PGF-dim}_RM = n$ implies that 
the inequality $\mbox{pd}_RK \leq n-1$ cannot be strict, i.e.\  
$\mbox{pd}_RK = n-1$. Since 
$K \in \overline{\tt Proj}(R) \subseteq {\tt GProj}(R)^{\perp} 
                              \subseteq {\tt PGF}(R)^{\perp}$,
where the latter inclusion is a consequence of the inclusion
${\tt PGF}(R) \subseteq {\tt GProj}(R)$, we conclude that $\pi$ 
is indeed a ${\tt PGF}(R)$-precover of $M$. \hfill $\Box$

\begin{Corollary}
If $M$ is a module with $\mbox{PGF-dim}_RM \leq 1$, then the 
following conditions are equivalent:

(i) $M \in {\tt PGF}(R)$,

(ii) $\mbox{Ext}^1_R(M,F)=0$ for any flat module $F$,

(iii) $\mbox{Ext}^1_R(M,P)=0$ for any projective module $P$.
\end{Corollary}
\vspace{-0.05in}
\noindent
{\em Proof.}
The implication (i)$\rightarrow$(ii) follows from [26, Theorem 4.4],
whereas the implication (ii)$\rightarrow$(iii) is obvious. In order 
to prove that (iii)$\rightarrow$(i), we use Proposition 3.1 and note 
that the hypothesis $\mbox{PGF-dim}_RM \leq 1$ implies the existence 
of a short exact sequence 
\[ 0 \longrightarrow P \longrightarrow G 
     \longrightarrow M \longrightarrow 0 , \]
where $P$ is projective and $G \in {\tt PGF}(R)$. By our assumption,
the group $\mbox{Ext}^1_R(M,P)$ is trivial and hence the exact 
sequence splits. It follows that $M$ is a direct summand of $G$. 
Since the class ${\tt PGF}(R)$ is closed under direct summands, 
we conclude that $M \in {\tt PGF}(R)$ as well. \hfill $\Box$

\begin{Corollary}
If $M \in \overline{\tt PGF}(R)$, then the following conditions are 
equivalent:

(i) $M \in {\tt PGF}(R)$,

(ii) $\mbox{Ext}^i_R(M,F)=0$ for any $i>0$ and any flat module $F$,

(iii) $\mbox{Ext}^i_R(M,P)=0$ for any $i>0$ and any projective module $P$.
\end{Corollary}
\vspace{-0.05in}
\noindent
{\em Proof.}
The implication (i)$\rightarrow$(ii) follows from [26, Corollary 4.5],
whereas the implication (ii)$\rightarrow$(iii) is obvious. In order to
prove that (iii)$\rightarrow$(i), we consider a PGF-resolution of $M$ 
of finite length
\[ 0 \longrightarrow G_n \longrightarrow \cdots \longrightarrow G_0 
     \longrightarrow M \longrightarrow 0 \]
and argue by induction on $n$. The case where $n=0$ is trivial. Assume 
that $n>0$ and let $K$ be the kernel of the map $G_0 \longrightarrow M$,
so that there is a short exact sequence
\[ 0 \longrightarrow K \longrightarrow G_0 \longrightarrow M 
     \longrightarrow 0 . \]
Since $G_0 \in {\tt PGF}(R)$, the group $\mbox{Ext}^i_R(G_0,P)$ is
trivial and hence $\mbox{Ext}^i_R(K,P) = \mbox{Ext}^{i+1}_R(M,P)$ is
also trivial for all $i>0$ and all projective modules $P$. The module 
$K$ admits a PGF-resolution of length $n-1$ and our induction hypothesis 
implies that $K \in {\tt PGF}(R)$. Therefore, it follows that
$\mbox{PGF-dim}_RM \leq 1$. Since $\mbox{Ext}^1_R(M,P) = 0$ for any 
projective module $P$, we finish the proof by invoking Corollary 3.2. 
\hfill $\Box$
\vspace{0.15in}
\newline
In view of Proposition 2.4(iii), the existence of a short exact 
sequence as in the statement of Proposition 3.1 is equivalent to 
the finiteness of the PGF-dimension of $M$. In fact, we may 
complement this assertion and prove the following result. Here, 
condition (iii) is analogous to [10, Lemma 2.17] (see also [25, 
Lemma 1.9]) and conditions (iv) and (v) are inspired by the Remark 
following [26, Theorem 4.11].

\begin{Theorem}
The following conditions are equivalent for a module $M$ and a 
non-negative integer $n$:

(i) $\mbox{PGF-dim}_RM = n$.

(ii) There exists a short exact sequence 
\[ 0 \longrightarrow K \longrightarrow G \longrightarrow M 
     \longrightarrow 0 , \]
where $G$ is a PGF-module and $\mbox{pd}_RK = n-1$. If $n=0$, this 
is understood to mean $K=0$. If $n=1$, we also require that the
exact sequence be non-split. 

(iii) There exists a short exact sequence 
\[ 0 \longrightarrow M \longrightarrow K \longrightarrow G 
     \longrightarrow 0 , \]
where $G$ is a PGF-module and $\mbox{pd}_RK = n$.

(iv) There exists a projective module $P$, such that the module 
$M' = M \oplus P$ fits into an exact sequence 
\[ 0 \longrightarrow G \longrightarrow M' \longrightarrow K 
     \longrightarrow 0 , \]
which remains exact after applying the functor 
$\mbox{Hom}_R(\_\!\_,Q)$ for any module 
$Q \in {\tt PGF}(R)^{\perp}$, where $G$ is a PGF-module and 
$\mbox{pd}_RK = n$.

(v) There exists a PGF-module $P$, such that the module 
$M' = M \oplus P$ fits into an exact sequence 
\[ 0 \longrightarrow G \longrightarrow M' \longrightarrow K 
     \longrightarrow 0 , \]
where $G$ is a PGF-module and $\mbox{pd}_RK = n$. If $n=1$,
we also require that the exact sequence remain exact after 
applying the functor $\mbox{Hom}_R(\_\!\_,Q)$ for any projective 
module $Q$. 
\end{Theorem}
\vspace{-0.05in}
\noindent
{\em Proof.}
(i)$\rightarrow$(ii): The existence of the short exact sequence 
follows from Proposition 3.1. If $n=1$, then the exact sequence 
cannot split. (Indeed, if the short exact sequence were split, 
then $M$ would be a direct summand of the PGF-module $G$ and 
hence $M$ would be itself a PGF-module; this is absurd, since 
$\mbox{PGF-dim}_RM=1$.)

(ii)$\rightarrow$(iii): Consider a short exact sequence as in (ii). 
Since $G \in {\tt PGF}(R)$, there exists a short exact sequence 
\[ 0 \longrightarrow G \longrightarrow P \longrightarrow G' 
     \longrightarrow 0 , \]
where $P$ is a projective module and $G' \in {\tt PGF}(R)$. By 
considering the pushout of that short exact sequence along the 
given epimorphism $G \longrightarrow M$, we obtain a commutative 
diagram with exact rows and columns
\[
\begin{array}{ccccccccc}
 & & & & 0 & & 0 & \\
 & & & & \downarrow & & \downarrow & \\
 0 & \longrightarrow & K & \longrightarrow & G 
   & \longrightarrow & M & \longrightarrow & 0 \\
 & & \parallel & & \downarrow & & \downarrow & \\
 0 & \longrightarrow & K & \longrightarrow & P 
   & \longrightarrow & K' & \longrightarrow & 0 \\
 & & & & \downarrow & & \downarrow & \\
 & & & & G' & = & G' & \\
 & & & & \downarrow & & \downarrow & \\
 & & & & 0 & & 0 & 
\end{array}
\]
We claim that the rightmost vertical exact sequence is of the 
required type. Indeed, if $n=0$, then $K=0$ and hence $K'=P$ is 
a projective module. If $n=1$, then $K$ is projective and the 
monomorphism $K \longrightarrow P$ is not split. (Indeed, if 
that monomorphism were split, then the monomorphism 
$K \longrightarrow G$ would be split as well, contradicting our 
assumption.) It follows that the module 
$K' = \mbox{coker} \left( K \longrightarrow P \right)$ is not
projective and hence $\mbox{pd}_RK'=1$. If $n \geq 2$, then 
$\mbox{pd}_RK=n-1>0$, so that 
$\mbox{Ext}^n_R(K',\_\!\_) = \mbox{Ext}^{n-1}_R(K,\_\!\_) \neq 0$ 
and $\mbox{Ext}^{n+1}_R(K',\_\!\_) = \mbox{Ext}^n_R(K,\_\!\_) = 0$;
it follows that $\mbox{pd}_RK'=n$.

(iii)$\rightarrow$(iv): Consider a short exact sequence as in 
(iii) and let 
\[ 0 \longrightarrow K' \longrightarrow P \longrightarrow K 
     \longrightarrow 0 \]
be a short exact sequence, where $P$ is a  projective module and 
$\mbox{pd}_RK' = n-1$. (If $n=0$, then $K$ is projective and we 
choose $P=K$ and $K'=0$.) By considering the pullback of that short 
exact sequence along the given monomorphism $M \longrightarrow K$, 
we obtain a commutative diagram with exact rows and columns
\[
\begin{array}{ccccccccc}
 & & 0 & & 0 & & & \\
 & & \downarrow & & \downarrow & & & \\
 & & K' & = & K' & \\
 & & \downarrow & & \downarrow & & & \\
 0 & \longrightarrow & G' & \longrightarrow & P 
   & \longrightarrow & G & \longrightarrow & 0 \\
 & & \downarrow & & \downarrow & & \parallel & \\
 0 & \longrightarrow & M & \longrightarrow & K 
   & \longrightarrow & G & \longrightarrow & 0 \\
 & & \downarrow & & \downarrow & & & \\
 & & 0 & & 0 & & &  
\end{array}
\]
Since the class ${\tt PGF}(R)$ is projectively resolving, the 
horizontal short exact sequence in the middle shows that $G'$
is a PGF-module. Then, the definition of the pullback and the 
surjectivity of the the map $P \longrightarrow K$ imply that 
there is a short exact sequence 
\[ 0 \longrightarrow G' \longrightarrow M \oplus P 
     \longrightarrow K \longrightarrow 0 . \]
In order to show that this short exact sequence has the required 
additional property, we note that for any module 
$Q \in {\tt PGF}(R)^{\perp}$ the two horizontal short exact 
sequences in the diagram above induce a commutative diagram 
of abelian groups with exact rows
\[
\begin{array}{ccccccccc}
 0 & \longrightarrow & \mbox{Hom}_R(G,Q) 
   & \longrightarrow & \mbox{Hom}_R(P,Q) 
   & \longrightarrow & \mbox{Hom}_R(G',Q) 
   & \longrightarrow & 0 \\
 & & \parallel & & \uparrow & & \uparrow & \\
 0 & \longrightarrow & \mbox{Hom}_R(G,Q) 
   & \longrightarrow & \mbox{Hom}_R(K,Q) 
   & \longrightarrow & \mbox{Hom}_R(M,Q) 
   & \longrightarrow & 0 
\end{array}
\]
It follows readily that there is an induced sequence of abelian 
groups 
\[ 0 \longrightarrow \mbox{Hom}_R(K,Q) 
     \longrightarrow \mbox{Hom}_R(M,Q) \oplus \mbox{Hom}_R(P,Q)  
     \longrightarrow \mbox{Hom}_R(G',Q) 
     \longrightarrow 0 , \]
as needed.

(iv)$\rightarrow$(v): This is immediate, since projective modules 
are contained in both classes ${\tt PGF}(R)$ and ${\tt PGF}(R)^{\perp}$.

(v)$\rightarrow$(i): Consider an exact sequence as in (v) and 
note that Proposition 2.3 implies that 
$\mbox{PGF-dim}_RM' = \mbox{PGF-dim}_RM$. Therefore, it suffices 
to prove that $\mbox{PGF-dim}_RM' = n$. Since $G$ is a PGF-module
and $\mbox{PGF-dim}_RK \leq \mbox{pd}_RK=n$, we may invoke Proposition 
2.4(i) and conclude that $\mbox{PGF-dim}_RM' \leq n$. It remains to 
show that the latter inequality cannot be strict. Indeed, let us 
assume that $n \geq 1$ and $\mbox{PGF-dim}_RM' \leq n-1$. 

If $n=1$, then $M'$ is a PGF-module and hence $\mbox{PGF-dim}_RK \leq 1$. 
Since the short exact sequence is assumed to remain exact after applying 
the functor $\mbox{Hom}_R(\_\!\_,Q)$ for any projective module $Q$ and 
$M' \in {\tt PGF}(R) \subseteq {^{\perp}{\tt Proj}(R)}$, it follows 
that the abelian group $\mbox{Ext}^1_R(K,Q)$ is trivial for any projective 
module $Q$. Then, Corollary 3.2 implies that $K \in {\tt PGF}(R)$; in 
particular, $K \in {\tt GProj}(R)$. As shown in [22, Proposition 2.27], 
any Gorenstein projective module of finite projective dimension is 
necessarily projective. We therefore conclude that the module $K$ is 
projective.\footnote{Alternatively, the projectivity of $K$ follows 
since $\overline{\tt Proj}(R) \subseteq {\tt PGF}(R)^{\perp}$ and 
${\tt PGF}(R) \cap {\tt PGF}(R)^{\perp} = {\tt Proj}(R)$; cf.\ [26].} 
This is absurd, since $\mbox{pd}_RK=1$.

We now consider the case where $n>1$. Since the PGF-module $G$ is
Gorenstein projective, the functor $\mbox{Ext}^{n-1}_R(G, \_\!\_)$ 
vanishes on projective modules. Since $\mbox{PGF-dim}_RM' \leq n-1$ 
and ${\tt PGF}(R) \subseteq {\tt GProj}(R)$, we also have 
$\mbox{Gpd}_RM' \leq n-1$. Therefore, [22, Theorem 2.20] implies that 
the functor $\mbox{Ext}^n_R(M', \_\!\_)$ vanishes on projective modules
as well. It follows that the functor $\mbox{Ext}^n_R(K,\_\!\_)$ vanishes
on projective modules. This contradicts our assumption that $\mbox{pd}_RK=n$; 
indeed, if we consider a projective resolution 
$P_* \longrightarrow K \longrightarrow 0$ of length $n$, then the 
monomorphism $P_n \longrightarrow P_{n-1}$ is not split and hence 
$\mbox{Ext}^n_R(K,P_n) \neq 0$. \hfill $\Box$
\vspace{0.15in}
\newline
{\bf Remarks 3.5.}
(i) In the case where $n=1$, it is necessary to impose some restrictions 
on the short exact sequences appearing in Theorem 3.4(ii),(v). Indeed, 
if $P$ is any non-zero projective module and $M \in {\tt PGF}(R)$, then 
the (split) short exact sequence 
\[ 0 \longrightarrow P \longrightarrow P \oplus M 
     \longrightarrow M \longrightarrow 0 \]
is of the type appearing in Theorem 3.4(ii), but 
$\mbox{PGF-dim}_RM = 0 \neq 1$. On the other hand, if $K$ is a module 
with $\mbox{pd}_RK=1$, then a projective resolution of $K$ provides an 
exact sequence 
\[ 0 \longrightarrow P_1 \longrightarrow P_0
     \longrightarrow K \longrightarrow 0 \]
of the type appearing in Theorem 3.4(v), but 
$\mbox{PGF-dim}_RP_0 = 0 \neq 1$.

(ii) It is clear from the proof of Theorem 3.4 that the analogues of 
conditions (iv) and (v) therein for Gorenstein projective modules are
equivalent to the analogues of conditions (i), (ii) and (iii) for such 
modules, thereby complementing the characterizations of the finiteness 
of the Gorenstein projective dimension given in [22, Theorem 2.10] and 
[10, Lemma 2.17].
\addtocounter{Lemma}{1}
\vspace{0.15in}
\newline
The next result is a characterization of modules of finite PGF-dimension,
that parallels the characterization of modules of finite Gorenstein 
projective dimension in [22, Theorem 2.20].

\begin{Proposition}
The following conditions are equivalent for a module $M$ of finite 
PGF-dimension and a non-negative integer $n$:

(i) $\mbox{PGF-dim}_RM \leq n$.

(ii) $\mbox{Ext}^i_R(M,F)=0$ for all $i>n$ and any flat module $F$.
     
(ii)' $\mbox{Ext}^i_R(M,P)=0$ for all $i>n$ and any projective module $P$.

(iii) $\mbox{Ext}^i_R(M,F)=0$ for all $i>n$ and any module $F$ of 
      finite flat dimension.

(iii)' $\mbox{Ext}^i_R(M,P)=0$ for all $i>n$ and any module $P$ of 
       finite projective dimension.
\end{Proposition}
\vspace{-0.05in}
\noindent
{\em Proof.}
(i)$\rightarrow$(ii): We consider a PGF-resolution of length $n$ 
\[ 0 \longrightarrow G_n \longrightarrow G_{n-1} 
     \longrightarrow \cdots \longrightarrow G_{0} 
     \longrightarrow M \longrightarrow 0 \]
and fix a flat module $F$. Since the functors $\mbox{Ext}^j_R(\_\!\_,F)$ 
vanish on the class of PGF-modules for all $j>0$ (cf.\ [26, Corollary 4.5]), 
we may deduce the desired vanishing by dimension shifting. 

(ii)$\rightarrow$(i): Let 
\[ 0 \longrightarrow K \longrightarrow G_{n-1} 
     \longrightarrow \cdots \longrightarrow G_0 
     \longrightarrow M \longrightarrow 0 \]
be an exact sequence, where $G_0, \ldots ,G_{n-1} \in {\tt PGF}(R)$. 
Since the modules $M,G_0, \ldots ,G_{n-1}$ are of finite PGF-dimension, 
an iterated application of Proposition 2.4(ii) shows that the module $K$ 
has finite PGF-dimension as well. On the other hand, our hypothesis and 
the dimension shifting argument employed in the proof of the implication 
(i)$\rightarrow$(ii) above show that the functors 
$\mbox{Ext}^i_R(K,\_\!\_)$ vanish on flat modules for all $i>0$. 
Invoking Corollary 3.3, we conclude that $K \in {\tt PGF}(R)$, 
as needed.

The implication (ii)$\rightarrow$(iii) follows by induction on the 
flat dimension of the module $F$, whereas the implication 
(iii)$\rightarrow$(ii) is immediate.

Finally, the implications 
(i)$\leftrightarrow$(ii)'$\leftrightarrow$(iii)' that involve projective 
modules can be proved by using exactly the same arguments as those used 
above for the implications that involve flat modules. \hfill $\Box$
\vspace{0.15in}
\newline
An immediate consequence of the characterization above is that the 
PGF-dimension is a refinement of the ordinary projective dimension,
whereas the Gorenstein projective dimension is a refinement of the 
PGF-dimension. 

\begin{Corollary}
Let $M$ be a module.

(i) If $\mbox{pd}_RM < \infty$, then $\mbox{PGF-dim}_RM = \mbox{pd}_RM$.

(ii) If $\mbox{PGF-dim}_RM < \infty$, then $\mbox{Gpd}_RM = \mbox{PGF-dim}_RM$.
\end{Corollary}
\vspace{-0.05in}
\noindent
{\em Proof.}
(i) If $\mbox{pd}_RM = n$, then the functors $\mbox{Ext}^i_R(M,\_\!\_)$
vanish for all $i>n$ and $\mbox{Ext}^i_R(M,P) \neq 0$ for a suitable 
projective module $P$. Since $\mbox{PGF-dim}_RM \leq n$, the equality 
$\mbox{PGF-dim}_RM=n$ follows from Proposition 3.6.

(ii) Since $\mbox{Gpd}_RM \leq \mbox{PGF-dim}_RM < \infty$, the 
equality $\mbox{Gpd}_RM = \mbox{PGF-dim}_RM$ follows from 
Proposition 3.6 and [22, Theorem 2.20]. \hfill $\Box$
\vspace{0.15in}
\newline
Since ${\tt PGF}(R) \subseteq {\tt GProj}(R)$, it follows from 
[22, Theorem 2.20] that $\mbox{Ext}^1_R(M,P) = 0$ whenever 
$M \in {\tt PGF}(R)$ and $P \in \overline{\tt Proj}(R)$; this is 
precisely the assertion of Proposition 3.6(iii)' in the case where 
$n=0$ therein. In fact, this vanishing provides a characterization 
of PGF-modules and modules of finite projective dimension, if we 
restrict to modules of finite PGF-dimension.

\begin{Proposition}
Let $N$ be a module of finite PGF-dimension. Then:

(i) $N \in {\tt PGF}(R)$ if and only if $\mbox{Ext}^1_R(N,P)=0$ 
    for any $P \in \overline{\tt Proj}(R)$. 
    
(ii) $N \in \overline{\tt Proj}(R)$ if and only if 
     $\mbox{Ext}^1_R(M,N)=0$ for any $M \in {\tt PGF}(R)$. 
\end{Proposition}
\vspace{-0.05in}
\noindent
{\em Proof.}
(i) As we noted above, the Ext-group is trivial if 
$N \in {\tt PGF}(R)$. Conversely, assume that $N$ is
a module of finite PGF-dimension contained in 
$^{\perp}\overline{\tt Proj}(R)$. Proposition 3.1 implies 
the existence of a short exact sequence 
\[ 0 \longrightarrow K \longrightarrow G \longrightarrow N
     \longrightarrow 0 , \]
where $G \in {\tt PGF}(R)$ and $K \in \overline{\tt Proj}(R)$. In
view of our assumption on $N$, this sequence splits and hence $N$
is a direct summand of the PGF-module $G$. Since the class 
${\tt PGF}(R)$ is closed under direct summands, we conclude that 
$N$ is a PGF-module.

(ii) As we noted above, the Ext-group is trivial if 
$N \in \overline{\tt Proj}(R)$. Conversely, assume that $N$ 
is module of finite PGF-dimension contained in 
${\tt PGF}(R)^{\perp}$. Then, Theorem 3.4(iii) implies the 
existence of a short exact sequence 
\[ 0 \longrightarrow N \longrightarrow K \longrightarrow G
     \longrightarrow 0 , \]
where $G \in {\tt PGF}(R)$ and $K \in \overline{\tt Proj}(R)$. In
view of our assumption on $N$, this sequence splits and hence $N$
is a direct summand of $K$. Then, 
$\mbox{pd}_RN \leq \mbox{pd}_RK < \infty$ and hence 
$N \in \overline{\tt Proj}(R)$, as needed. \hfill $\Box$
\vspace{0.15in}
\newline
We now examine the special case of Gorenstein flat modules and 
show that the values of their PGF-dimension are controlled by the 
values of the projective dimension of flat modules. We let 
$\mbox{splf} \, R$ be the supremum of the projective lengths 
(dimensions) of flat modules.

\begin{Proposition}
We have an equality 
$\sup \{ \mbox{PGF-dim}_RM : M \in {\tt GFlat}(R) \} = \mbox{splf} \, R$.
In particular, ${\tt Flat}(R) \subseteq \overline{\tt Proj}(R)$
if and only if ${\tt GFlat}(R) \subseteq \overline{\tt PGF}(R)$.
\end{Proposition}
\vspace{-0.05in}
\noindent
{\em Proof.}
Let $s = \sup \{ \mbox{PGF-dim}_RM : M \in {\tt GFlat}(R) \}$. 
If $M$ is any flat module, then Proposition 2.5(i) implies that 
$\mbox{pd}_RM = \mbox{PGF-dim}_RM \leq s$. It follows that 
$\mbox{splf} \, R \leq s$. In order to prove the reverse inequality,
it suffices to assume that $\mbox{splf} \, R < \infty$, so that any 
flat module has finite projective dimension. If $M$ is any Gorenstein 
flat module, then [26, Theorem 4.11] implies that there exists a short 
exact sequence 
\[ 0 \longrightarrow M \longrightarrow F \longrightarrow G
     \longrightarrow 0 , \]
where $F$ is flat and $G \in {\tt PGF}(R)$. Since $F$ has finite 
projective dimension, Theorem 3.4(iii) implies that 
$\mbox{PGF-dim}_RM = \mbox{pd}_RF \leq \mbox{splf} \, R$. We 
conclude that $s \leq \mbox{splf} \, R$, as needed. 

Considering the projective dimension of direct sums of flat modules,
it is easily seen that ${\tt Flat}(R) \subseteq \overline{\tt Proj}(R)$
if and only if $\mbox{splf} \, R < \infty$. In the same way, we may 
consider the PGF-dimension of direct sums of Gorenstein flat modules 
(cf.\ Proposition 2.3) and conclude that 
${\tt GFlat}(R) \subseteq \overline{\tt PGF}(R)$ if and only if 
$s < \infty$. Therefore, the final statement in the Proposition 
follows from the equality $s = \mbox{splf} \,R$. \hfill $\Box$
\vspace{0.15in}
\newline
We may complement the characterization of the finiteness of 
PGF-dimension given in Theorem 3.4, in the case of a Gorenstein
flat module $M$, by requiring that the module $K$ that appears 
in assertions (ii), (iii), (iv) and (v) therein be also flat. To 
that end, we note that any Gorenstein flat module of finite projective 
dimension is necessarily flat. Indeed, such a module must also have 
finite flat dimension and its flatness follows then from [6, $\S 2$]; 
see also [15, Remark 1.5].

\begin{Proposition}
The following conditions are equivalent for a Gorenstein flat 
module $M$ and a non-negative integer $n$:

(i) $\mbox{PGF-dim}_RM = n$.

(ii) There exists a short exact sequence 
\[ 0 \longrightarrow K \longrightarrow G \longrightarrow M 
     \longrightarrow 0 , \]
where $G$ is a PGF-module and $K$ is a flat module with 
$\mbox{pd}_RK = n-1$. If $n=0$, this is understood to mean $K=0$. 
If $n=1$, we also require that the exact sequence be non-split. 

(iii) There exists a short exact sequence 
\[ 0 \longrightarrow M \longrightarrow K \longrightarrow G 
     \longrightarrow 0 , \]
where $G$ is a PGF-module and $K$ is a flat module with 
$\mbox{pd}_RK = n$.

(iv) There exists a projective module $P$, such that the module 
$M' = M \oplus P$ fits into an exact sequence 
\[ 0 \longrightarrow G \longrightarrow M' \longrightarrow K 
     \longrightarrow 0 , \]
which remains exact after applying the functor 
$\mbox{Hom}_R(\_\!\_,Q)$ for any module 
$Q \in {\tt PGF}(R)^{\perp}$, where $G$ is a PGF-module and 
$K$ is a flat module with $\mbox{pd}_RK = n$.

(v) There exists a PGF-module $P$, such that the module 
$M' = M \oplus P$ fits into an exact sequence 
\[ 0 \longrightarrow G \longrightarrow M' \longrightarrow K 
     \longrightarrow 0 , \]
where $G$ is a PGF-module and $K$ is a flat module with 
$\mbox{pd}_RK = n$. If $n=1$, we also require that the exact 
sequence remain exact after applying the functor 
$\mbox{Hom}_R(\_\!\_,Q)$ for any projective module $Q$. 
\end{Proposition}
\vspace{-0.05in}
\noindent
{\em Proof.}
We proceed as in the proof of Theorem 3.4, showing that
(i)$\rightarrow$(ii)$\rightarrow$(iii)$\rightarrow$(iv)$\rightarrow$(v)$\rightarrow$(i). 
Since $M$ is Gorenstein flat, ${\tt PGF}(R) \subseteq {\tt GFlat}(R)$ 
and the class of Gorenstein flat modules is projectively resolving
(cf.\ [26, Corollary 4.12]), the module $K$ appearing in (ii) and (iii) 
is a Gorenstein flat module of finite projective dimension; as noted
above, this forces $K$ to be flat. We also note that the argument 
in the proof of the implication (iii)$\rightarrow$(iv) in Theorem 
3.4 provides a short exact sequence as in (iv) with $K$ being {\em 
the same} module $K$ that appears in (iii). \hfill $\Box$

\section{Hovey triples on $\overline{\tt PGF}(R)$ and ${\tt GFlat}(R)$}

\noindent
We shall now relate the results obtained in the previous section to 
the theory of exact model structures and describe a hereditary Hovey 
triple in the exact category $\overline{\tt PGF}(R)$ of modules of 
finite PGF-dimension, which is such that the homotopy category of 
the associated exact model structure is equivalent as a triangulated 
category to the stable category of PGF-modules. We shall also describe 
the stable category of PGF-modules, up to triangulated equivalence, as 
the homotopy category of the exact model structure associated with a 
similar Hovey triple in the exact category ${\tt GFlat}(R)$ of 
Gorenstein flat modules.

It is easily seen that ${\tt PGF}(R)$ is an exact Frobenius category 
with projective-injective objects given by the projective modules. 
The proof of the latter claim is essentially identical to the proof 
of the corresponding claim for the class of Gorenstein projective 
modules, which can be found for instance in [13, Proposition 2.2]. 

The category $\overline{\tt PGF}(R)$ of modules of finite PGF-dimension 
is an extension closed subcategory of the abelian category of all modules 
(cf.\ Proposition 2.4(i)), which is also closed under direct summands
(cf.\ Proposition 2.3). Therefore, $\overline{\tt PGF}(R)$ is an idempotent 
complete exact additive category [8]. The following result is an analogue 
of [13, Theorem 3.7]. The idea is that in order to realize the stable 
category of PGF-modules as the homotopy category of a Quillen model 
structure, it suffices to work on the subcategory $\overline{\tt PGF}(R)$ 
of modules of finite PGF-dimension. We note that the class 
$\overline{\tt Proj}(R)$ of modules of finite projective dimension 
is closed under direct summands and has the 2-out-of-3 property for 
short exact sequences. 

\begin{Theorem}
The triple 
$\left( {\tt PGF}(R) , \overline{\tt Proj}(R) , 
 \overline{\tt PGF}(R) \right)$ 
is a hereditary Hovey triple in the idempotent complete 
exact category $\overline{\tt PGF}(R)$. The homotopy category 
of the associated exact model structure is equivalent, as a 
triangulated category, to the stable category of PGF-modules.
\end{Theorem}
\vspace{-0.05in}
\noindent
{\em Proof.} 
We need to prove that the pairs 
\[ \left( {\tt PGF}(R) , 
   \overline{\tt Proj}(R) \cap \overline{\tt PGF}(R) \right)
   \;\;\; \mbox{and} \;\;\;
   \left( {\tt PGF}(R) \cap \overline{\tt Proj}(R) , 
   \overline{\tt PGF}(R) \right) \]
are complete and hereditary cotorsion pairs in the exact category 
$\overline{\tt PGF}(R)$. Since any PGF-module is Gorenstein projective,
we conclude that  
\[ {\tt PGF}(R) \cap \overline{\tt Proj}(R) \subseteq 
   {\tt GProj}(R) \cap \overline{\tt Proj}(R) = {\tt Proj}(R) , \]
where the latter equality follows from [22, Proposition 2.27]. On 
the other hand, projective modules are contained in both classes 
${\tt PGF}(R)$ and $\overline{\tt Proj}(R)$ and hence 
${\tt PGF}(R) \cap \overline{\tt Proj}(R) = {\tt Proj}(R)$.\footnote{Alternatively, 
the equality ${\tt PGF}(R) \cap \overline{\tt Proj}(R) = {\tt Proj}(R)$
follows since 
${\tt Proj}(R) \subseteq \overline{\tt Proj}(R) \subseteq 
 {\tt PGF}(R)^{\perp}$ 
and ${\tt PGF}(R) \cap {\tt PGF}(R)^{\perp} = {\tt Proj}(R)$; cf.\ [26].}
Thus, the two pairs displayed above become
\[ \left( {\tt PGF}(R) , \overline{\tt Proj}(R) \right)
   \;\;\; \mbox{and} \;\;\;
   \left( {\tt Proj}(R) , \overline{\tt PGF}(R) \right) . \]

We begin by considering the pair 
$\left( {\tt PGF}(R) , \overline{\tt Proj}(R) \right)$ and note that
Proposition 3.8 states precisely that this is indeed a cotorsion pair
in $\overline{\tt PGF}(R)$. Theorem 3.4 provides the approximations 
referring to completeness, whereas Proposition 3.6(iii)', applied to 
the case where $n=0$, shows that the cotorsion pair is hereditary.

We now consider the pair 
$\left( {\tt Proj}(R) , \overline{\tt PGF}(R) \right)$ and note 
that $\overline{\tt PGF}(R)$ is obviously the right orthogonal of 
${\tt Proj}(R)$ within $\overline{\tt PGF}(R)$. In order to prove 
that ${\tt Proj}(R)$ is the left orthogonal of
$\overline{\tt PGF}(R)$ within $\overline{\tt PGF}(R)$, we let 
$M$ be a module of finite PGF-dimension which is also contained in 
$^{\perp}\overline{\tt PGF}(R)$ and consider a short exact sequence
\[ 0 \longrightarrow M' \longrightarrow P \longrightarrow M
     \longrightarrow 0 , \]
where $P$ is projective. Then, Proposition 2.4(ii) implies that $M'$ 
has also finite PGF-dimension and hence $\mbox{Ext}^1_R(M,M') = 0$.
In particular, the exact sequence above splits. It follows that $M$ 
is a direct summand of $P$ and hence $M$ is projective. The cotorsion 
pair $\left( {\tt Proj}(R) , \overline{\tt PGF}(R) \right)$ in 
$\overline{\tt PGF}(R)$ is hereditary (since all higher Ext's with 
a projective first argument vanish) and complete (since the class
$\overline{\tt PGF}(R)$ is projectively resolving). 

The rest of the statement follows from [18, Proposition 4.4 and 
Corollary 4.8]. \hfill $\Box$
\vspace{0.15in}
\newline
The category ${\tt GFlat}(R)$ of Gorenstein flat modules is also
closed under extensions and direct summands; this follows from 
[26, Corollary 4.12]. Hence, ${\tt GFlat}(R)$ is an idempotent 
complete exact category as well. As shown in [26, Theorem 4.4], the 
group $\mbox{Ext}^1_R(M,F)$ is trivial whenever $M$ is a PGF-module 
and $F$ is flat. This vanishing actually provides a characterization 
of PGF-modules and flat modules, if we restrict to Gorenstein flat 
modules. 

\begin{Proposition}
Let $N$ be a Gorenstein flat module. Then:

(i) $N \in {\tt PGF}(R)$ if and only if $\mbox{Ext}^1_R(N,F)=0$ 
    for any flat module $F$.
        
(ii) $N$ is flat if and only if $\mbox{Ext}^1_R(M,N)=0$ for any 
     $M \in {\tt PGF}(R)$.
\end{Proposition}
\vspace{-0.05in}
\noindent
{\em Proof.}
(i) As we noted above, the Ext-group is trivial if $N$ is a PGF-module. 
Conversely, assume that $N$ is a Gorenstein flat module contained 
in $^{\perp}{\tt Flat}(R)$. Then, there exists a short exact sequence 
\[ 0 \longrightarrow F \longrightarrow G \longrightarrow N
     \longrightarrow 0 , \]
where $G$ is a PGF-module and $F$ is flat; cf.\ [26, Theorem 4.11(2)].
In view of our assumption on $N$, this short sequence splits and hence 
$N$ is a direct summand of the PGF-module $G$. Since the class 
${\tt PGF}(R)$ is closed under direct summands, we conclude that 
$N \in {\tt PGF}(R)$.

(ii) As we noted above, the Ext-group is trivial if $N$ is flat. 
Conversely, assume that $N$ is a Gorenstein flat module contained 
in ${\tt PGF}(R)^{\perp}$. Then, there exists a short exact sequence 
\[ 0 \longrightarrow N \longrightarrow F \longrightarrow G
     \longrightarrow 0 , \]
where $G$ is a PGF-module and $F$ is flat; cf.\ [26, Theorem 4.11(4)]. 
In view of our assumption on $N$, this short sequence splits and hence 
$N$ is a direct summand of the flat module $F$. Therefore, $N$ is flat. 
\hfill $\Box$
\vspace{0.15in}
\newline
We note that the class ${\tt Flat}(R)$ of flat modules is closed 
under direct summands and has the 2-out-of-3 property within the 
class of Gorenstein flat modules. Of course, ${\tt Flat}(R)$ is 
closed under extensions and kernels of epimorphisms. Moreover, if 
the cokernel of a monomorphism between flat modules is Gorenstein 
flat, then that cokernel is necessarily flat.\footnote{We have pointed 
out in the discussion preceding Proposition 3.10 that any Gorenstein 
flat module of finite flat dimension is necessarily flat.} The proof 
of the following result is very similar to the proof of Theorem 4.1. 

\begin{Theorem}
The triple 
$\left( {\tt PGF}(R) , {\tt Flat}(R) , {\tt GFlat}(R) \right)$ is 
a hereditary Hovey triple in the idempotent complete exact category 
${\tt GFlat}(R)$. The homotopy category of the associated exact model 
structure is equivalent, as a triangulated category, to the stable 
category of PGF-modules.
\end{Theorem}
\vspace{-0.05in}
\noindent
{\em Proof.} 
We need to prove that the pairs 
\[ \left( {\tt PGF}(R) , {\tt Flat}(R) \cap {\tt GFlat}(R) 
   \right)
   \;\;\; \mbox{and} \;\;\;
   \left( {\tt PGF}(R) \cap {\tt Flat}(R) , {\tt GFlat}(R)
   \right) \]
are complete and hereditary cotorsion pairs in the exact category 
${\tt GFlat}(R)$. Since 
${\tt PGF}(R) \cap {\tt Flat}(R) = {\tt Proj}(R)$, the two pairs 
displayed above become
\[ \left( {\tt PGF}(R) , {\tt Flat}(R) \right)
   \;\;\; \mbox{and} \;\;\;
   \left( {\tt Proj}(R) , {\tt GFlat}(R) \right) . \]

We begin by considering the pair 
$\left( {\tt PGF}(R) , {\tt Flat}(R) \right)$ and note that Proposition 
4.2 states precisely that this is indeed a cotorsion pair in the exact
category ${\tt GFlat}(R)$. Completeness of the cotorsion pair follows 
from the exact sequences in [26, Theorem 4.11(2),(4)], whereas Proposition 
3.6(ii), applied to the case where $n=0$, shows that the cotorsion pair 
is hereditary.

We now consider the pair 
$\left( {\tt Proj}(R) , {\tt GFlat}(R) \right)$ and note that 
${\tt GFlat}(R)$ is obviously the right orthogonal of ${\tt Proj}(R)$ 
within ${\tt GFlat}(R)$. In order to prove that ${\tt Proj}(R)$ is the 
left orthogonal of ${\tt GFlat}(R)$ within ${\tt GFlat}(R)$, we let 
$M$ be a Gorenstein flat module which is also contained in 
$^{\perp}{\tt GFlat}(R)$ and consider a short exact sequence
\[ 0 \longrightarrow M' \longrightarrow P \longrightarrow M
     \longrightarrow 0 , \]
where $P$ is projective. Since the class ${\tt GFlat}(R)$ is projectively 
resolving (cf.\ [26, Corollary 4.12]), we deduce that $M'$ is also Gorenstein 
flat. Therefore, $\mbox{Ext}^1_R(M,M') = 0$ and the exact sequence above 
splits. It follows that $M$ is a direct summand of $P$ and hence $M$ is 
projective. The cotorsion pair 
$\left( {\tt Proj}(R) , {\tt GFlat}(R) \right)$ in ${\tt GFlat}(R)$ is 
hereditary (since all higher Ext's with a projective first argument 
vanish) and complete (since the class of Gorenstein flat modules is 
projectively resolving). 

The final statement follows from [18, Proposition 4.4 and Corollary 4.8]. 
\hfill $\Box$
\vspace{0.15in}
\newline
{\bf Remark 4.4.}
Another model for the stable category of PGF-modules can be obtained 
from the Hovey triple 
$({\tt PGF}(R),{\tt PGF}(R)^{\perp},\mbox{$R$-Mod})$ on the category
$R$-Mod of all modules; cf.\ [26, Theorem 4.9] and [19, Proposition 37]. 
A possible advantage of the Hovey triples presented in this section 
is that the classes of modules that are involved herein admit a more 
manageable description.
\addtocounter{Lemma}{1}

\section{The finiteness of the PGF global dimension}

\noindent
In this section, we characterize those rings over which all modules 
have finite PGF-dimension, in terms of classical homological invariants. 
As a consequence of this description, we generalize a result by Jensen 
[24] (on commutative Noetherian rings) and another result by Gedrich and 
Gruenberg [17] (on group rings of groups over a commutative Noetherian 
coefficient ring).

We define the (left) PGF global dimension $\mbox{PGF-gl.dim} \, R$ 
of the ring $R$, by letting 
\[ \mbox{PGF-gl.dim} \, R = \sup \{ \mbox{PGF-dim}_{R}M : 
   \mbox{$M$ a left $R$-module} \} . \]
Using the characterization of the finiteness of the Gorenstein global 
dimension and the Goresntein weak global dimension, we may characterize 
the finiteness of $\mbox{PGF-gl.dim} \, R$, as follows:

\begin{Theorem}
The following conditions are equivalent for a ring $R$:

(i) $\mbox{PGF-gl.dim} \, R < \infty$,

(ii) $\mbox{PGF-dim}_RM < \infty$ for any module $M$,

(iii) $\mbox{spli} \, R = \mbox{silp} \, R < \infty$ and 
      $\mbox{sfli} \, R = \mbox{sfli} \, R^{op} < \infty$,

(iv) $\mbox{spli} \, R < \infty$ and 
     $\mbox{sfli} \, R^{op} < \infty$.
\newline
If these conditions are satisfied, then
$\mbox{PGF-gl.dim} \, R = \mbox{spli} \, R = 
 \mbox{silp} \, R ( = \mbox{Ggl.dim} \, R )$.
\end{Theorem}
\vspace{-0.05in}
\noindent
{\em Proof.}
It is clear that (i)$\rightarrow$(ii), whereas the implication 
(ii)$\rightarrow$(i) is an immediate consequence of Proposition 2.3.

(ii)$\rightarrow$(iii): Since ${\tt PGF}(R)$ is contained in both 
classes ${\tt GProj}(R)$ and ${\tt GFlat}(R)$, our hypothesis implies 
that any module $M$ has both finite Gorenstein projective dimension 
and finite Gorenstein flat dimension. (In fact, both $\mbox{Gpd}_RM$ 
and $\mbox{Gfd}_RM$ are bounded by $\mbox{PGF-dim}_RM < \infty$.) Then, 
assertion (iii) follows from the characterization of the finiteness of 
the Gorenstein global dimension and the Gorenstein weak global dimension 
of $R$; cf.\ $\S 1$.II.

(iii)$\rightarrow$(iv): This is straightforward.

(iv)$\rightarrow$(ii): Assume that $\mbox{spli} \, R = n < \infty$ and 
fix a module $M$. Then, the construction by Gedrich and Gruenberg in 
[17, $\S 4$] provides us with an acyclic complex of projective modules 
\begin{equation}
 \cdots \longrightarrow P_{n+1} \longrightarrow P_{n}
 \longrightarrow Q_{n-1} \longrightarrow Q_{n-2} \longrightarrow
 \cdots , 
\end{equation}
which coincides in degrees $\geq n$ with a projective resolution
\[ \cdots \longrightarrow P_{n+1} \longrightarrow P_{n}
   \longrightarrow P_{n-1} \longrightarrow \cdots
   \longrightarrow P_{0} \longrightarrow M \longrightarrow 0 \]
of $M$. Since the acyclic complex (3) consists of projective (and hence 
flat) modules, it remains acyclic by applying the functor 
$L \otimes_R \_\!\_$ for any right module $L$ of finite flat dimension;
this follows easily by induction on the flat dimension of $L$. In 
particular, our assumption about the finiteness of $\mbox{sfli} \, R^{op}$ 
implies that the complex (3) remains acyclic by applying the functor 
$I \otimes_R \_\!\_$ for any injective right module $I$. Therefore, the 
module $K = \mbox{coker} \left( P_{n+1} \longrightarrow P_n \right)$ is 
a PGF-module. Then, the exact sequence 
\[ 0 \longrightarrow K \longrightarrow P_{n-1} \longrightarrow 
     \cdots \longrightarrow P_{0} \longrightarrow M \longrightarrow 0 \]
shows that $\mbox{PGF-dim}_RM \leq n < \infty$, as needed. 

The final claim in the statement of the Theorem is an immediate 
consequence of Corollary 3.7(ii), which implies that 
$\mbox{PGF-gl.dim} \, R = \mbox{Ggl.dim} \, R$, if
$\mbox{PGF-gl.dim} \, R$ is finite. \hfill $\Box$
\vspace{0.15in}
\newline
As an immediate consequence of the equivalence between assertions (iii) 
and (iv) in Theorem 5.1 above, we obtain the following result.

\begin{Corollary}
Let $R$ be a ring, such that both invariants $\mbox{spli} \, R$ and
$\mbox{sfli} \, R^{op}$ are finite. Then, 
$\mbox{silp} \, R = \mbox{spli} \, R$. \hfill $\Box$
\end{Corollary}

\noindent
We may obtain a left-right symmetric assertion, as follows.

\begin{Proposition}
Let $R$ be a ring, such that both invariants $\mbox{spli} \, R$ and
$\mbox{spli} \, R^{op}$ are finite. Then, we have
$\mbox{silp} \, R = \mbox{spli} \, R$ and 
$\mbox{silp} \, R^{op} = \mbox{spli} \, R^{op}$.
\end{Proposition}
\vspace{-0.05in}
\noindent
{\em Proof.}
Since projective (left or right) modules are flat, we have
\[ \mbox{sfli} \, R \leq \mbox{spli} \, R < \infty
   \;\;\; \mbox{and} \;\;\;
   \mbox{sfli} \, R^{op} \leq \mbox{spli} \, R^{op} < \infty . \]
Then, the result follows by applying Corollary 5.2 for the ring 
$R$ and its opposite $R^{op}$. \hfill $\Box$

\begin{Corollary}
If $R$ is a ring which is isomorphic with its opposite $R^{op}$, 
then $\mbox{silp} \, R \leq \mbox{spli} \, R$ with equality if 
$\mbox{spli} \, R < \infty$.
\end{Corollary}
\vspace{-0.05in}
\noindent
{\em Proof.}
The inequality is obvious if $\mbox{spli} \, R = \infty$ and 
hence it suffices to consider the case where 
$\mbox{spli} \, R < \infty$. Then, 
$\mbox{spli} \, R^{op} = \mbox{spli} \, R$ is also finite and 
we may invoke Proposition 5.3. \hfill $\Box$
\vspace{0.15in}
\newline
We recall that a ring $R$ is called left (resp.\ right) 
$\aleph_0$-Noetherian if any left (resp.\ right) ideal of $R$ is 
countably generated. For example, countable rings and countably 
generated algebras over fields are both left and right 
$\aleph_0$-Noetherian.
\vspace{0.15in}
\newline
{\bf Remarks 5.5.}
(i) Let $k$ be a commutative ring, $G$ a group and $R=kG$ the 
associated group algebra. Then, $R$ is isomorphic with its 
opposite $R^{op}$ and hence Corollary 5.4 implies that 
$\mbox{silp} \, R \leq \mbox{spli} \, R$. In the special case 
where the coefficient ring $k$ is Noetherian of finite self-injective
dimension, this inequality was proved by Gedrich and Gruenberg 
in [17, Theorem 2.4], using the Hopficity of the group algebra 
$R$.

(ii) Let $k$ be a commutative $\aleph_0$-Noetherian ring, 
$G$ a group and $R=kG$ the associated group algebra. Then, 
we may invoke [14, Proposition 4.3] and conclude that the 
inequality in (i) above is actually an equality, i.e.\ 
$\mbox{silp} \, R = \mbox{spli} \, R$. In this way, we
extend the main result of [14] from the case of commutative 
Noetherian rings of finite self-injective dimension to any 
commutative $\aleph_0$-Noetherian ring of coefficients. 
\addtocounter{Lemma}{1}

\begin{Proposition}
If $R$ is a ring which is both left and right $\aleph_0$-Noetherian, 
then the following conditions are equivalent:

(i) The invariants $\mbox{spli} \, R$ and
    $\mbox{spli} \, R^{op}$ are finite. 

(ii) The invariants $\mbox{silp} \, R$ and
     $\mbox{silp} \, R^{op}$ are finite.
\newline
If these conditions are satisfied, then 
$\mbox{silp} \, R = \mbox{spli} \, R < \infty$ and 
$\mbox{silp} \, R^{op} = \mbox{spli} \, R^{op} < \infty$.
\end{Proposition}
\vspace{-0.05in}
\noindent
{\em Proof.}
The implication (i)$\rightarrow$(ii) follows from Proposition 
5.3, whereas the implication (ii)$\rightarrow$(i) is proved in 
[14, Theorem 3.6] using the $\aleph_0$-Noetherian hypothesis.
\hfill $\Box$

\begin{Corollary}
Let $R$ be a ring which is isomorphic with its opposite $R^{op}$.
If $R$ is left (and hence right) $\aleph_0$-Noetherian, then 
$\mbox{silp} \, R = \mbox{spli} \, R$. \hfill $\Box$
\end{Corollary}

\begin{Corollary}
If $R$ is a commutative $\aleph_0$-Noetherian ring, then 
$\mbox{silp} \, R = \mbox{spli} \, R$. \hfill $\Box$
\end{Corollary}

\noindent
{\bf Remarks 5.9.}
(i) In the special case where $R$ is a commutative Noetherian ring, 
the equality in Corollary 5.8 was proved by Jensen in [24, 5.9].

(ii) The analogous result of Proposition 5.6 for left and right 
coherent rings appears in [3, Theorem 3.3].
\addtocounter{Lemma}{1}
\vspace{0.1in}
\newline
{\em Acknowledgments.}
We thank the anonymous referee who read carefully the manuscript
and offered helpful comments.

\vspace{0.15in}

\noindent
{\small {\sc Institute of Algebra and Number Theory, 
             University of Stuttgart, Pfaffenwaldring 57, 
             70569 Stuttgart, Germany}}

\noindent
{\em E-mail address:} {\tt gdalezios@math.uoa.gr} 

\vspace{0.05in}
\noindent
and

\vspace{0.05in}
\noindent
{\small {\sc Department of Mathematics,
             University of Athens,
             Athens 15784,
             Greece}}

\noindent
{\em E-mail address:} {\tt emmanoui@math.uoa.gr}

\end{document}